\input amstex 
\documentstyle{amsppt}
\input bull-ppt
\keyedby{bull281/lic}
\define\otr{\operatorname{Tr}}
\define\vl{\operatorname{Vol}}
\define\pa{\partial}
\define\trv{|\!|\!|}
\topmatter
\cvol{26}
\cvolyear{1992}
\cmonth{April}
\cyear{1992}
\cvolno{2}
\cpgs{299-303}

\title Spectral Theory and Representations\\ 
of Nilpotent Groups\endtitle
\author P. Levy-Bruhl, A. Mohamed, and J. Nourrigat 
\endauthor
\shortauthor{P. Levy-Bruhl, A. Mohamed, and J. Nourrigat}
\shorttitle{Spectral theory}
\address {\rm(P. Levy-Bruhl and J. Nourrigat)} 
D\'epartement de Mathematiques, Universit\'e de Reims, B. 
P. 347,
51062 Reims CEDEX, France \endaddress
\address {\rm(A. Mohamed)} D\'epartment de Mathematiques, 
Universit\'e de Nantes,
2 rue de la Houssini\`ere, 44072 Nantes CEDEX, 
France\endaddress
\date May 15, 1991 and, in revised form, September 24, 
1991\enddate
\thanks The results have been presented at the 1991 
meeting on
``Harmonic analysis and representation theory'' (R. S. 
Howe. org.) (held January
27 to February 2),
in Oberwolfach, Germany\endthanks
\abstract We give an estimate of the number $N(\lambda)$ of 
eigenvalues $<\lambda$ for the image under an irreducible 
representation
of the ``sublaplacian'' on a stratified nilpotent Lie 
algebra. We
also give an estimate for the trace of the heat-kernel 
associated
with this operator. The estimates are formulated in term 
of geometrical
objects related to the representation under consideration. 
An important
particular case is the Schr\"odinger equation with 
polynomial 
electrical and magnetical fields.\endabstract
\endtopmatter

\document
\heading 1. Introduction\endheading
We first consider the Schr\"odinger operator with 
polynomial electrical
and magnetical fields
$$
P=\sum_{j=1}^n (D_j-A_j(x))^2 +V(x),
$$
where $A_j$ and $V$ are real polynomials on $\Bbb R^n$, of 
degree
$\le r$, and with $V\ge 0$. We define $B_{jk}=\pa A_j/\pa 
x_k-
\pa A_k/\pa x_j$ and assume that there is no rotation of the
coordinates axis making all the $B_{jk}$ and $V$ independent
of one of the $x_i$'s. This case has been investigated in 
\cite9:
define
$$
\gather
M(x)=\sum_{|\alpha|\le r} |\pa^\alpha V(x)|^{1/(|\alpha|+2)}
+\sum_{\alpha,j,k} |\pa^\alpha B_{jk}(x)|^{1/(|\alpha|+
2)},\\
M(x,\xi)=|\xi|+M(x)\.
\endgather
$$

Let $N(\lambda)$ be the number of eigenvalues of 
$P<\lambda$ and 
$N_0(\lambda)$ the volume in $\Bbb R^{2n}$ of the set of 
points
$(x,\xi)$ such that $M(x,\xi)^2\le\lambda$. Then there is a 
$C>1$, independent of $\lambda$ such that $C^{-1} N_0(C^{-1}
\lambda)\le N(\lambda)\le CN_0(C\lambda)$ for all 
$\lambda>0$.
A more precise equivalent in particular cases is given in 
\cite{13}.
This result was proved before, in a different formulation, 
by
Fefferman, when the $A_j$'s are zero \cite 1. If $V$ is 
the square
of a polynomial, $P$ is of the form $\pi(-\Delta)$ where 
$\Delta$
is a sublaplacian on a stratified $r$-step nilpotent Lie 
group
and $\pi$ an irreducible representation \cite{2, 9}. In the
general case, $(V\ge 0)$, our following
results\ needs a very small
change to be applied. This example is one of the 
motivations for the
present generalization.

\heading 2. Statement of the results \endheading
Let $\germ G$ be a stratified $r$-steps nilpotent Lie 
algebra. 
In other words, we assume that $\germ G$ can be written as 
a direct
sum of subspaces $\germ G_j$ $(1\le j\le r)$ such that
$$\gather
[\germ G_j,\germ G_k]\subset \germ G_{j+k}\quad\text{if}\ 
j+k\le r,\\
[\germ G,\germ G_k]=0\quad\text{if}\ j+k>r,
\endgather
$$
and such that $\germ G$ is generated as Lie algebra by the 
subspace
$\germ G_1$. We choose a basis $X_1,\dots X_p$ of $\germ 
G_1$, and
we shall be interested in spectral properties of the image 
under
irreducible representations of the sublaplacian 
$\Delta=-\sum_{j=1}^p
 X_j^2$. 

Let $\pi$ be an irreducible, unitary, nontrivial 
representation of the
connected simply connected group $\exp\germ G$ associated 
to $\germ G$.
It is well known that the operator $\pi(-\Delta)$, defined 
on the space
$\scr S_\pi$ of $C^\infty$ vectors of the representation 
$\pi$, has a
unique selfadjoint extension, still denoted 
$\pi(-\Delta)$, whose 
spectrum is a sequence $(\lambda_j)$ $(j\in\Bbb N)$ of 
positive 
eigenvalues, such that $\lambda_j\le\lambda_{j+1}$ and 
$\lambda_j\to
+\infty$. Let us denote by $N(\lambda)$ the number of 
eigenvalues
$\lambda_j\le\lambda$. Our first goal is to give an 
estimate for $N
(\lambda)$ in terms of geometrical objects associated to 
the 
representation $\pi$. 

We denote by $\delta_t$ $(t>0)$ the natural dilations of 
$\germ G$;
that is, the linear maps defined by $\delta_t(X)=t^jX$ if 
$X\in\germ G_j$.
Let also $\germ G^*$ be the dual of $\germ G$, 
$\delta_t^*$ the
transpose of $\delta_t$, and $\trv\ \trv$ a homogeneous 
norm on $\germ G^*$,
i.e., a nonnegative continuous, subadditive function on 
$\germ G^*$,
only vanishing at the origin and satisfying 
$\trv\delta_t^*l\trv=
t\trv l\trv$ if $t>0$ and $l\in\germ G^*$. 

By the Kirillov theory, the representation $\pi$ is 
associated with an
orbit $O_\pi$ of the coadjoint representation of the group 
$G=\exp
\germ G$ in $\germ G^*$. The $G$ invariant measures on 
$O_\pi$ are
proportional, and we denote by $\mu$ the canonical one; 
that is, 
with the correct normalization for the character formula 
\cite{4,
12}. For each $\lambda>0$, we set
$$
N_0(\lambda)=\mu(l\in O_\pi,\ \trv l\trv^2\le \lambda)\.
$$

We can now state the estimate for $N(\lambda)$.

\thm{Theorem 1} There exists a constant $C>1$, independent 
of $\lambda$,
and of the representation $\pi$ such that
$$
C^{-1} N_0(C^{-1}\lambda) \le N(\lambda)\le 
CN_0(C\lambda),\quad
\text{for all}\ \lambda>0\.
$$
\ethm

Manchon gives in \cite7 the proof of a conjecture of 
Karasev-Maslov
\cite3. His result does not apply to our operator, but 
rather to
$-\sum_j Y_j^2$, where $(Y_j)$ is a basis of $\germ G$ as 
a vector
space. Such an operator is also selfadjoint with compact 
resolvent,
and $N(\lambda)$ is given by a Weyl type formula. Such a 
formula 
does not make sense in general in our case, the integral 
involved 
in it being nonconvergent.

It is also well known, since $\pi$ is irreducible, that 
for $t>0$,
the operator $\exp(-t\pi(\Delta))$ is trace-class. We can 
estimate
its trace
$$
Z(t)=\otr (\exp(-t\pi(\Delta)))
$$
by the following result, using the function
$$
Z_0(t)=\int_{\scr O_\pi} \exp(-t\trv l\trv^2)\,d\mu(l)\.
$$

\thm{Theorem 2} There exists a constant $C>1$, independent 
of $t>0$,
and of the representation $\pi$ such that
$$
C^{-1} Z_0(Ct) \le Z(t)\le CZ_0(C^{-1} t)\.
$$
\ethm

As noticed in the introduction, our results can be used 
for the
Schr\"odinger operator with polynomial electrical and 
magnetical fields.

\heading 3. Sketch of the proof\endheading

\rem{Remark} We give here a sketch of the proof of Theorem 
1. The reader may
consult \cite6 for the details.\endrem

We construct Hilbert spaces of sequences controlling the 
norm of the
Sobolev spaces associated with $\pi(-\Delta)$ and use the 
minimax
formula.

For this construction, we need a suitable set of functions 
in $L^2
(\Bbb R^n)$. The properties of these functions have been 
suggested
both by the ideas of Perelomov \cite{11} (coherent states) 
and of
Meyer \cite8 (wavelets). 

We realize the representation $\pi$ in the following form: 
for 
$X\in\germ G$, $\pi(X)$ is a differential operator on 
$\Bbb R^n$,
$$\aligned
\pi(X)=& A_1(X)\pa/\pa x_1+A_2(x_1,X)\pa/\pa x_2 \\
&+\cdots+A_n(x_1,\dots, x_{n-1},X)\pa/\pa x_n+iB(x,X),
\endaligned\tag "$(*)$"
$$
where the $A_j$ and $B$ are real linear forms in $X$ and 
polynomials in
$x\in\Bbb R^n$, $A_1$ is independent of $x$, and $A_j$ 
only depends on
$(x_1,\dots, x_{j-1})$. For every finite sequence 
$I=(i_1,\dots, i_m)$
of positive integers smaller than $p$, $\pi(X_I)$ is the 
iterated
commutator
$$
\pi(X_I)=(\roman{ad}\,\pi (X_{i_1}))\cdots (\roman{ad}\, 
\pi(X_{i_{m-1}}))
\pi (X_{i_m})\.
$$
Let $\pi(X)(x,\xi)$ be the complete symbol of the operator 
$\pi(X)$.
The orbit $\scr O_\pi$ is the set of linear forms 
$X\to-i\pi (X)(x,\xi)$
for $(x,\xi)$ in $\Bbb R^{2n}$, and the canonical measure 
$\mu$ is then
$dxd\xi$. We can also choose the homogeneous norm
$$
\trv l\trv=\sum |l(X_i)|^{1/|I|}\.
$$
The functions $M_\pi(x,\xi)$ and $M_\pi(x)$ are defined by
$$
M_\pi(x,\xi)=\sum_{|I|\le r} 
|\pi(X_I)(x,\xi)|^{1/|I|}\quad\text{and}\quad
M_\pi(x)=\inf_\xi M_\pi(x,\xi)\.
$$

For the representation $\pi$ under consideration, 
$M_\pi(x)>0$. We
prove the theorems using this \RM{``}concrete\RM{''} 
realization of $\pi$. 
We are not able to construct an orthonormal basis of $L^2$ 
so that
the Sobolev spaces associated to our operator are $l^2$ 
weighted
spaces, but the following result is a substitute and 
allows the use
of minimax formulas. Let $p(x,\xi)=x$, $(x,\xi)\in\Bbb 
R^n$. For
every finite sequence $(\alpha_1,\dots, \alpha_m)$ of 
positive integers
smaller than $p$, we define 
$\pi(X^\alpha)=\pi(X_{\alpha_1})\cdots
\pi(X_{\alpha_m})$, and for
$$
u\in C_0^\infty (\Bbb R^n)\: \|u\|_{m,\pi}
=\lf(\sum_{|\alpha|=m} \|\pi(X^\alpha) u\|^2\rt)^{1/2}\.
$$

\thm{Theorem 3} 
\RM{``}Coherent states.\RM{''} Let $\pi$ be an induced
representation of $\germ G$ realized in the form $(*)$ with
$M_\pi(x)>0$ for all $x\in\Bbb R^n$. \RM(\<$\pi$ is not 
necessarily
irreducible.\RM) For every $a>0$ small enough, there 
exists $(\psi_j)$
and $(\Psi_j)$ $(j\in\Bbb N)$ in $C_0^\infty(\Bbb R^n)$, a 
sequence
$\Omega_j$ of compact sets in $\Bbb R^{2n}$, a point 
$(x_j,\xi_j)\in
\Omega_j$, and $C>0$ with: 

\RM1. $u=\sum_j (u,\Psi_j)\psi_j$ for all $u$ in $L^2(\Bbb 
R^n)$.

\RM2. $(\supp \Psi_j)\subset(\supp)\psi_j\subset 
p(\Omega_j)$ 
for all $j\in\Bbb N$. 

\RM3. Every finite subset $A$ of $\Bbb N$ contains a part 
$B$
such that $\#B\ge (\#A)/C$ and such that for every sequence
$(\lambda_j)$ of complex numbers
$$
\sum_{j\in B}|\lambda_j|^2 \le C\lf\|\sum_{j\in B} \lambda_j
\psi_j\rt\|^2\.
$$

\RM4. For every integer $m\ge 0$, there exists $C_m>0$ 
with the 
properties\,\RM: If $u\in C_0^\infty(\Bbb R^n)$ is of the 
form
$u=\sum_j\lambda_j\psi_j$ $(\lambda_j\in\Bbb C$ not 
necessarily
given by $(u,\Psi_j)$\<\RM) then $\|u\|_{m,\pi}^2\le 
C_m\sum_j |\lambda
_j|^2 M_\pi(x_j,\xi_j)^{2m}$. If in addition 
$\lambda_j=(u,\Psi_j)$ 
then we have $\sum_j|(u,\Psi_j)|^2 
M_\pi(x_j,\xi_j)^{2m}\le C_m
\|u\|_{m,\pi}^2$. 

\RM5. $\Bbb R^{2n}\subset\bigcup_j (\Omega_j)$ and every 
$(x,\xi)$
is contained in at most $C$ of the $\Omega_j$\RM's. 
Further, $1/C\le
\vl(\Omega_j)\le C$. 

\RM6. For $(x,\xi)$ in $\Omega_j$, we have the estimate 
$(1/C)M_\pi
(x_j,\xi_j)\le M_\pi(x,\xi)\le CM_\pi(x_j,\xi_j)$.

\RM7. $C$ and the constants $C_m$ only depend on $a$ and 
$\germ G$
not on $\pi$.\ethm

The construction of the coherent states is first performed 
locally:
we associate to each $x\in\Bbb R^n$ a representation 
$\sigma_x$,
equivalent to $\pi$ such that $\sigma_x$ has the form 
$(*)$ with 
$M_{\sigma_x}(0)=M_\pi(x)$ and the coefficients of the 
polynomials
in the expression of $\sigma_x(\delta_t^{-1}(X))$ where 
$t=M_\pi(x)$
are bounded, independently of $x$ and $\pi$. Using Fourier 
series
and induction on the dimension $n$, ``local coherent 
states'' for
$\sigma_x$ are constructed. Certain properties in Theorem 
3 are
only satisfied if $u$ is supported in the ball of radius 
$a$,
centered at 0. The coherent states are then obtained by 
means of 
a partition of unity, modeled on the symplectic 
diffeomorphisms
related to the intertwining operators between $\pi$ and 
$\sigma_x$. 

In the proof, the following theorem is used (the notation 
is as
before).

\thm{Proposition 4} There exists $C>0$ such that for every 
representation $\pi$ of the form $(*)$ and all $u\in 
S(\Bbb R^n)$
$$
\sum_{|\alpha|\le m} \|M_\pi(x)^{m-|\alpha|} \pi(X^\alpha) 
u\|^2
\le C\|u\|_{m,\pi}^2\.
$$
\ethm

In this form, this theorem requires $\germ G$ to be 
stratified,
however, the proof of Lemma 4.5 in \cite{10} works without 
this
hypothesis. The function $M$ is then more complicated, 
reflecting
the structure of $\germ G$. It is possible to formulate 
our present
results in this set-up. The same method also applies to more
general operators than the model presented here; we need 
only a
control of $(\pi(P)u,u)$ in terms of the $\|\ \|_{m,\pi}$ 
norm.

\enddocument